\documentclass[12pt]{amsart}
\usepackage{latexsym,amssymb,amscd,eqname}

\def\B'c{{\mathcal{B'}}}
\def\U'c{{\mathcal{U'}}}

\def\opn#1#2{\def#1{\operatorname{#2}}} 
\opn\chara{char}
\opn\length{\ell}
\opn\projdim{proj\,dim}
\opn\injdim{inj\,dim}
\opn\ini{in}
\opn\rank{rank}
\opn\depth{depth}
\opn\height{ht}
\opn\embdim{emb\,dim}
\opn\codim{codim}

\opn\Tr{Tr}
\opn\bigrank{big\,rank}
\opn\superheight{superheight}\opn\lcm{lcm}
\opn\trdeg{tr\,deg}%
\opn\reg{reg}
\opn\lreg{lreg}
\opn\set{set}
\opn\supp{Supp}
\opn\shad{Shad}
\opn\div{div}
\opn\Div{Div}
\opn\cl{cl}
\opn\Cl{Cl}
\opn\Spec{Spec}
\opn\Supp{Supp}
\opn\supp{supp}
\opn\Sing{Sing}
\opn\Ass{Ass}
\opn\Ann{Ann}
\opn\Rad{Rad}
\opn\Soc{Soc}
\opn\Ker{Ker}
\opn\Coker{Coker}
\opn\Im{Im}
\opn\Hom{Hom}
\opn\Tor{Tor}
\opn\Ext{Ext}
\opn\End{End}
\opn\Aut{Aut}
\opn\id{id}

\opn\nat{nat}
\opn\GL{GL}
\opn\SL{SL}
\opn\mod{mod}
\opn\ord{ord}
\opn\aff{aff}
\opn\con{conv}
\opn\relint{relint}
\opn\st{st}
\opn\lk{lk}
\opn\cn{cn}
\opn\core{core}
\opn\vol{vol}
\opn\gr{gr}

\def\pot#1#2{#1[\kern-0.28ex[#2]\kern-0.28ex]}

\opn\dirlim{\underrightarrow{\lim}}
\opn\invlim{\underleftarrow{\lim}}
\def\pnt{{\raise0.5mm\hbox{\large\bf.}}}

\def\Implies{\ifmmode\Longrightarrow \else
     \unskip${}\Longrightarrow{}$\ignorespaces\fi}
\def\implies{\ifmmode\Rightarrow \else
     \unskip${}\Rightarrow{}$\ignorespaces\fi}
\def\iff{\ifmmode\Longleftrightarrow \else
     \unskip${}\Longleftrightarrow{}$\ignorespaces\fi}

\let\:=\colon
\newtheorem{Theorem}{Theorem}[section]
\newtheorem{Lemma}[Theorem]{Lemma}
\newtheorem{Corollary}[Theorem]{Corollary}
\newtheorem{Proposition}[Theorem]{Proposition}

\newtheorem{Example}[Theorem]{Example}

\newtheorem{Definition}[Theorem]{Definition}

\let\epsilon=\varepsilon
\let\phi=\varphi
\let\kappa=\varkappa
\numberwithin{equation}{section}

\title{Properties of lexsegment ideals }
\author{Viviana Ene \and Anda Olteanu \and Loredana Sorrenti}
\thanks{\footnotesize  The first author was supported by the grant CEX 05-D11-11/2005.
The second author was supported by the CNCSIS grant TD 507/2007.
}

\address{Faculty of Mathematics and Computer Science, Ovidius University, Bd.\ Mamaia 124,
 900527 Constanta, Romania,} \email{vivian@univ-ovidius.ro} 
\address{Faculty of Mathematics and Computer Science, Ovidius University, Bd.\ Mamaia 124,
 900527 Constanta, Romania,} \email{olteanuandageorgiana@gmail.com} 
\address{DIMET University of Reggio Calabria, Faculty of Engineering, via Graziella (Feo di Vito), 89100 Reggio Calabria, Italy} \email{loredana.sorrenti@unirc.it. }

\begin{document}

\maketitle

\begin{abstract} We show that any lexsegment ideal with linear resolution has linear quotients with respect to a suitable ordering of its minimal monomial generators. For completely lexsegment ideals with linear resolution we show that the decomposition function is regular. For arbitrary lexsegment ideals we compute the depth and the dimension. As application we characterize the Cohen-Macaulay lexsegment ideals. \\



\end{abstract}

\section*{Introduction}

Let $S=k[x_1,\ldots, x_n]$ be the polynomial ring in $n$ variables over a field $k$. We order lexicographically the monomials of $S$ such that $x_1>x_2>\ldots> x_n$. Let $d\geq 2$ be an integer and $\mathcal{M}_d$ the set of monomials of degree $d$. For two monomials $u,v\in\mathcal{M}_d$, with $u\geq_{lex}v$, the set $$\mathcal L(u,v)=\{w\in\mathcal{M}_d\ |\ u\geq_{lex}w\geq_{lex}v\}$$ is called a lexsegment. A lexsegment ideal in $S$ is a monomial ideal of $S$ which is generated by a lexsegment. Lexsegment ideals have been introduced by Hulett and Martin \cite{HM}. Arbitrary lexsegment ideals have been studied by A. Aramova, E. De Negri, and J. Herzog in \cite{ADH} and \cite{DH}. They characterized the lexsegment ideals which have linear resolutions.

In this paper we show that any lexsegment ideal with linear resolution has linear quotients with respect to a suitable order of the generators.

Let $I\subset S$ be a monomial ideal and $G(I)$ its minimal monomial set of generators. $I$ has linear quotients if there exists an ordering $u_1,\ldots,u_m$ of the elements of $G(I)$ such that for all $2\leq j\leq m$, the colon ideals $(u_1,\ldots,u_{j-1}):u_j$ are generated by a subset of $\{x_1,\ldots,x_n\}$.

Lexsegment ideals which have linear quotients with respect to the lexicographical order of the generators have been characterized by the third author in \cite{S}.

In  Section \ref{Section1} we show that any completely lexsegment ideal with linear resolution has linear quotients with respect to the following order of the generators. Given two monomials of degree $d$ in $S$, $w=x_1^{\alpha_1}\ldots x_n^{\alpha_n}$ and $w'=x_1^{\beta_1}\ldots x_n^{\beta_n},$   we set $w\prec w'$ if $\alpha_1<\beta_1$ or $\alpha_1=\beta_1$ and $w>_{lex}w'$.

Let $u,v\in\mathcal{M}_d$ which define the completely lexsegment ideal $I=(\mathcal{L}(u,v))$ with linear resolution. If $\mathcal{L}(u,v)=\{w_1,\ldots, w_r\}$, where $w_1\prec w_2\prec\ldots\prec w_r$, we show that $I$ has linear quotients with respect to this ordering of the generators. The non-completely lexsegment ideal will be separately studied in Section $2$.

For the completely lexsegment ideals with linear resolution it will turn out that their decomposition function with respect to the ordering $\prec$  is regular . Therefore, one may apply the procedure developed in \cite{HT} to get the explicit resolutions for this class of ideals.

In the last section of our paper we study the depth and the dimension of lexsegment ideals. Our results show that one may compute these invariants just looking at the ends of the lexsegment. As an application, we characterize the Cohen-Macaulay lexsegment ideals.\\

We acknowledge the support provided by the Computer Algebra Systems \textsc{CoCoA} \cite{Co} and \textsc{Singular} \cite{GPS} for the extensive experiments which helped us to obtain some of the results of this work. 

\section{Completely lexsegment ideals with linear resolutions}\label{Section1}

In the theory of Hilbert functions or in extremal combinatorics usually one considers initial lexsegment ideals, that is ideals generated by an initial lexsegment $\mathcal{L}^i(v)=\{w\in\mathcal{M}_d\ |\ w\geq_{lex} v\}$.   Initial lexsegment ideals are stable in the sense of Eliahou and Kervaire (\cite{EH}, \cite{AH}) and they have linear quotients with respect to lexicographical order  \cite[Proposition 2.1]{S}.

One may also define the final lexsegment $\mathcal{L}^f(u)=\{w\in\mathcal{M}_d\ |\ u\geq_{lex}w\}$. Final lexsegment ideals are generated by final lexsegments. They are also stable in the sense of Eliahou and Kervaire  with respect to $x_n>x_{n-1}>\ldots>x_1$. Therefore they have linear quotients.

Throughout this paper we use the following notations. If $m=x_1^{\alpha_1}\ldots x_n^{\alpha_n}$ is a monomial of $S$, we denote by $\nu_i(m)$ the exponent of the variable $x_i$ in $m$, that is $\nu_i(m)=\alpha_i$, $i=1,\ldots,n$. Also, we will denote $\max(m)=\max\{i\mid x_i|m\}.$

Hulett and Martin call a lexsegment $L$  \it{completely lexsegment} \rm if all the iterated shadows of $L$ are again lexsegments. We recall that the shadow of a set $T$ of monomials is the set $\shad(T)=\{vx_i\mid v\in T,\ 1\leq i\leq n\}$. The $i$-th shadow is recursively defined as $\shad^i(T)=\shad(\shad^{i-1}(T))$. The initial lexsegments have the property that their shadow is again an initial lexsegment, a fact which is not true for arbitrary lexsegments. An ideal spanned by a completely lexsegment is called a \it completely lexsegment ideal\rm. All the completely lexsegment ideals with linear resolution are determined in \cite{ADH}:

\begin{Theorem}\cite{ADH}\label{completelylex} Let $u=x_1^{a_1}\ldots x_n^{a_n}$, $v=x_1^{b_1}\ldots x_n^{b_n}$ be monomials of degree $d$ with $u\geq_{lex}v$, and let $I=(\mathcal{L}(u,v))$ be a completely lexsegment ideal. Then $I$ has linear resolution if and only if one of the following conditions holds:
\begin{itemize}
	\item [(a)] $u=x_1^ax_2^{d-a},\ v=x_1^ax_n^{d-a}$ for some $a,\ 0< a\leq d;$
	\item [(b)] $b_1< a_1-1;$
	\item [(c)] $b_1 = a_1-1$ and for the largest $w <_{lex} v, w$ monomial of degree $d,$ one has $x_1 w/x_{\max(w)}\leq_{lex} u.$
\end{itemize}
\end{Theorem}

\begin{Theorem}\label{colex}
Let $u=x_1^{a_1}\ldots x_n^{a_n},$ with  $a_1>0,$ and $ v=x_1^{b_1}\ldots x_n^{b_n}$ be monomials of degree $d$ with $u\geq_{lex} v,$ and let $I=({\mathcal L}(u,v))$ be a completely lexsegment ideal. Then $I$ has linear resolution if and only if $I$ has linear quotients.
\end{Theorem}

\begin{proof}
We have to prove that if $I$ has linear resolution then $I$ has linear quotients, since the other implication is known \cite{H}. By Theorem \ref{completelylex}, since $I$ has linear resolution, one of the conditions (a), (b), (c) holds.

We define on the set of the monomials of  degree $d$ from $S$ the following total order: for $$w=x_1^{\alpha_1}\ldots x_n^{\alpha_n},\ w^{\prime}=x_1^{\beta_1}\ldots x_n^{\beta_n},$$ we set $$w\prec w^{\prime}\text{ if }\alpha_1 < \beta_1 \text{ or }\alpha_1=\beta_1 \text{ and } w>_{lex} w^{\prime}.$$ Let $$\mathcal{L}(u,v)=\{w_1,\ldots,w_r\}, \text{ where }w_1\prec w_2\prec\ldots\prec w_r.$$ We will prove that $I=(\mathcal{L}(u,v))$ has linear quotients with respect to this ordering of the generators. 

Assume that $u,v$ satisfy the condition (a) and $a<d$ (the case $a=d$ is trivial). Then $I$ is isomorphic as $S$--module with the ideal generated by the final lexsegment ${\mathcal L}^f(x_2^{d-a}) \subset S$ and the ordering $\prec$ of its minimal generators coincides with the lexicographical ordering $>_{lex}.$ The ideal $(\mathcal{L}^f(x_2^{d-a}))\cap k[x_2,\ldots,x_n]$ is the initial ideal in $k[x_2,\ldots,x_n]$ defined by $x_n^{d-a},$ which has linear quotients with respect to $>_{lex}.$ Hence $I$ has linear quotients with respect to $\prec$ since it is obvious that the extension in the ring $k[x_1,\ldots, x_n]$ of a monomial  ideal with linear quotients in $k[x_2,\ldots,x_n]$ has linear quotients too.

Next we assume that $u,v$ satisfy the condition (b) or (c).

By definition, $I$ has linear quotients with respect to the monomial generators $w_1,\ldots,w_r$ if the colon ideals $(w_1,\ldots,w_{i-1}):w_i$ are generated by variables for all $i\geq 2,$ that is for all $j < i$ there exists an integer $1\leq k <i$ and an integer $l\in [n]$ such that $w_k/\gcd(w_k,w_i)=x_l \text{ and }x_l\text{ divides } w_j/\gcd(w_j,w_i).$ 

In other words, for any $w_j\prec w_i, w_j,w_i\in \mathcal{L}(u,v),$ we have to find a monomial $w^{\prime}\in \mathcal{L}(u,v)$ such that 
\begin{eqnarray}
w^{\prime}\prec w_i,\ \frac{w^{\prime}}{\gcd(w^{\prime},w_i)}=x_l, \text{ for some }l\in [n]  \text{, and }x_l\text{ divides } \frac{w_j}{\gcd(w_j,w_i)}. \eqname{$*$}\label{*}
\end{eqnarray}

Let us fix  $w_i=x_1^{\alpha_{1}}\ldots x_n^{\alpha_{n}}$ and $w_j=x_1^{\beta_{1}}\ldots x_n^{\beta_{n}}, \ w_i,w_j\in \mathcal{L}(u,v),$ such that $w_j\prec w_i.$ By the definition of the  ordering $\prec,$ we must have $$\beta_{1} < \alpha_{1} \text { or }\beta_{1} = \alpha_{1} \text{ and } w_j>_{lex} w_i.$$

\underline{\textit{Case $1$}}: Let $\beta_1<\alpha_1$. One may find an integer $l$, $2\leq l\leq n$, such that $\alpha_s\geq \beta_s$ for all $s<l$ and $\alpha_l<\beta_l$ since, otherwise, $\deg(w_i)>\deg(w_j)=d$ which is impossible. We obviously have $\max(w_j)\geq l$. If $l\geq \max(w_i),$ one may take $\bar{w}=x_lw_i/x_1$ which satisfies the condition (\ref{*}) since the inequalities $\bar{w}\prec w_i,$ $\bar{w}\leq_{lex}w_i\leq_{lex} u$ hold, and we will show that $\bar{w}\geq_{lex}w_j$. This will imply that $\bar{w}\geq_{lex}v$, hence $\bar{w}\in\mathcal{L}(u,v)$.

The inequality $\bar{w}\geq_{lex}w_j$ is obviously fulfilled if $\alpha_1-1>\beta_1$ or if $\alpha_1-1=\beta_1$ and at least one of the inequalities $\alpha_s\geq\beta_s$ for $2\leq s<l$, is strict. If $\alpha_1-1=\beta_1$ and $\alpha_s=\beta_s$ for all $s<l$, comparing the degrees of $w_i$ and $w_j$ it results $d=\alpha_1+\ldots+\alpha_l=\beta_1+1+\beta_2+\ldots+\beta_{l-1}+\alpha_l<(\beta_1+1)+\beta_2+\ldots+\beta_l.$ It follows that $d\geq\beta_1+\beta_2+\ldots+\beta_l>d-1,$ that is $\beta_1+\beta_2+\ldots+\beta_l=d.$ This implies that $l=\max(w_j)$ and $\beta_l=\alpha_l+1$, that is $\bar{w}=x_lw_i/x_1=x_1^{\alpha_1-1}x_2^{\alpha_2}\ldots x_l^{\alpha_l+1}=x_1^{\beta_1}\ldots x_l^{\beta_l}=w_j.$ 

From now on, in the Case $1$, we may assume that $l<\max(w_i)$. We will show that at least one of the following monomials: $$w'=\frac{x_lw_i}{x_{\max(w_i)}},\ w''=\frac{x_lw_i}{x_1}$$ belongs to $\mathcal{L}(u,v)$. It is clear that both monomials are strictly less than $w_i$ with respect to the ordering $\prec.$ Therefore one of the monomials $w',\ w''$ will satisfy the condition (\ref{*}).

The following inequalities are fulfilled: $$w'>_{lex}w_i\geq_{lex}v,\ \text{and}$$ $$w''<_{lex}w_i\leq_{lex}u.$$

Let us assume, by contradiction, that $w'>_{lex}u$ and $w''<_{lex}v$. Comparing the exponents of the variable $x_1$, we obtain $a_1-1\leq\alpha_1-1\leq b_1$. Since the ideal generated by $\mathcal{L}(u,v)$ has linear resolution, we must have $b_1=a_1-1$. Let $z$ be the largest monomial of degree $d$ such that $z<_{lex}v$. Then, by our assumption on $w''$, we also have the inequality $w''\leq_{lex}z$.

Now we need the following

\begin{Lemma}\label{lema1} Let $m=x_1^{a_1}\ldots x_n^{a_n},\ m'=x_1^{b_1}\ldots x_n^{b_n}$ be two monomials of degree $d$. If $m\leq_{lex}m'$ then $m/x_{\max(m)}\leq_{lex}m'/x_{\max(m')}$.
\end{Lemma}
\begin{proof} Let $m <_{lex}m'$. Then there exists $s\geq1$ such that $a_1=b_1,\ldots, a_{s-1}=b_{s-1}$ and $a_s<b_s$. It is clear that $\max(m')\geq s$. Comparing the degrees of $m$ and $m'$  we get  $\max(m)>s.$ 

If $\max(m)>s$ and $\max(m')>s$, the required inequality is obvious.

Let $\max(m)>s$ and $\max(m')=s$. Let us suppose, by contradiction, that $m/x_{\max(m)}>_{lex}m'/x_{\max(m')}=m'/x_s.$ This implies that $a_s\geq b_s-1$, and, since $a_s<b_s$, we get $a_s=b_s-1$.
Looking at the degree of $m'$ we obtain $d=b_1+b_2+\ldots+b_s=a_1+a_2+\ldots+ a_{s-1}+a_s+1,$ that is $a_1+\ldots+a_s=d-1$. It follows that $a_{\max(m)}=1$ and $m/x_{\max(m)}=x_1^{a_1}\ldots x_s^{a_s}=x_1^{b_1}\ldots x_{s-1}^{b_{s-1}}x_s^{b_s-1}=m'/x_{\max(m')}$, contradiction.
\end{proof}
Going back to the proof of our theorem, 
we apply the above lemma for the monomials $w''$ and $z$ and we obtain  $w''/x_{\max(w'')}\leq_{lex}z/x_{\max(z)},$ which implies that $x_1w''/x_{\max(w'')}\leq_{lex}x_1z/x_{\max(z)}.$
By using  condition (c) in the Theorem \ref{completelylex} it follows that $x_1w''/x_{\max(w'')}\leq_{lex}u.$ On the other hand, $x_1w''/x_{\max(w'')}=x_1x_lw_i/(x_1x_{\max(w_i)})$ $=x_lw_i/x_{\max(w_i)}=w'.$ Therefore, it results $w'\leq_{lex}u$, which contradicts our assumption on $w'$.

Consequently, we have $w'\leq_{lex}u$ or $w''\geq_{lex}v$, which proves that at least one of the monomials $w',\ w''$ belongs to $\mathcal{L}(u,v)$.

\underline{\textit{Case $2$:}} Let $\beta_1=\alpha_1$ and $w_j>_{lex}w_i$. Then there exists $l$, $2\leq l\leq n$, such that $\alpha_s=\beta_s$, for all $s<l$ and $\alpha_l<\beta_l$. If $\max(w_i)\leq l$, then, looking at the degrees of $w_i$ and $w_j$, we get $d=\alpha_1+\alpha_2+\ldots+\alpha_l<\beta_1+\beta_2+\ldots+\beta_l,$ contradiction. Therefore, $l<\max(w_i)$. We proceed in a similar way as in the previous case. Namely, exactly as in the Case $1,$ it results that at least one of the following two monomials $w'=x_lw_i/x_{\max(w_i)},\ w''=x_lw_i/x_1$ belongs to $\mathcal{L}(u,v)$. It is clear that both monomials are strictly less than $w_i$ with respect to the order $\prec$.  
\end{proof}

\begin{Example}\rm 
Let $S=k[x_1,x_2,x_3]$. We consider the monomials: $u=x_1x_2x_3$ and $v=x_2x_3^2$, $u>_{lex} v,$ and let $I$ be the monomial ideal generated by $\mathcal{L}(u,v)$. The minimal system of generators of the ideal $I$ is
	$$G(I)=\mathcal{L}(u,v)=\{x_1x_2x_3,\ x_1x_3^2,\ x_2^3,\ x_2^2x_3,\ x_2x_3^2\}.$$

Since $I$ verifies the condition (c) in Theorem \ref{completelylex}, it follows that $I$ is a completely lexsegment ideal with linear resolution. We denote the monomials from $G(I)$ as follows: $u_1=x_1x_2x_3,\ u_2=x_1x_3^2,\ u_3=x_2^3,\ u_4=x_2^2x_3,\ u_5=x_2x_3^2$, so $u_1>_{lex}u_2>_{lex}\ldots>_{lex}u_5.$ The colon ideal $(u_1,u_2):u_3=(x_1x_3)$ is not generated by a subset of $\{x_1,x_2,x_3\}$. This shows that $I$ is not with linear quotients with respect to lexicographical order.

We consider now the  order $\prec$ and check by direct computation that $I$ has linear quotients.  We label the monomials from $G(I)$ as follows: $u_1=x_2^3,\ u_2=x_2^2x_3,\ u_3=x_2x_3^2,\ u_4=x_1x_2x_3,\ u_5=x_1x_3^2$, so $u_1\prec u_2\prec\ldots\prec u_5$. Then $(u_1):u_2=(x_2),\ (u_1,u_2):u_3=(x_2),\ (u_1,u_2,u_3):u_4=(x_2,x_3),$ $(u_1,u_2,u_3,u_4):u_5=(x_2)$. 
\end{Example}

We further study the decomposition function of a completely lexsegment ideal with linear resolution. The decomposition function of a monomial ideal was introduced by J. Herzog and Y. Takayama in \cite{HT}. 

We recall the following notation. If $I\subset S$ is a monomial ideal with linear quotients with respect to the ordering $u_1,\ldots,u_m$ of its minimal generators, then we denote $$\set(u_j)=\{k\in[n]\ |\ x_k\in(u_1,\ldots,u_{j-1}):u_j\}$$
for $j=1,\ldots,m$.

\begin{Definition}\cite{HT} \rm Let $I\subset S$ be a monomial ideal with linear quotients with respect to the sequence of minimal monomial generators $u_1,\ldots, u_m$ and set $I_j=(u_1,\ldots, u_j)$, for $j=1,\ldots,m$. Let $\mathcal{M}(I)$ be the set of all monomials in $I$. The map $g:M(I)\rightarrow G(I)$ defined as: $g(u)=u_j$, where $j$ is the smallest number such that $u\in I_j$, is called \it the decomposition function \rm of $I$.
\end{Definition}

We say that the decomposition function $g:M(I)\rightarrow G(I)$ is \it regular \rm if $\set(g(x_su))$$\subseteq\set(u)$ for all $s\in\set(u)$ and $u\in G(I)$.

We show in the sequel that completely lexsegment ideals which have linear quotients with respect to $\prec$ have also regular decomposition functions.

In order to do this, we need some preparatory notations and results.

For an arbitrary lexsegment $\mathcal{L}(u,v)$ with the elements ordered by $\prec$, we denote by $I_{\prec w}$, the ideal generated by all the monomials $z\in\mathcal{L}(u,v)$ with $z\prec w$. $I_{\preceq}w$ will be the ideal generated by all the monomials $z\in\mathcal{L}(u,v)$ with $z\preceq w$.

\begin{Lemma}\label{1nset} Let $I=(\mathcal{L}(u,v))$ be a  lexsegment ideal which has linear quotients with respect to the order $\prec$ of the generators. Then, for any $w\in\mathcal{L}(u,v)$, $1\notin\set(w)$.
\end{Lemma}

\begin{proof} Let us assume that $1\in\set(w)$, that is $x_1w\in I_{\prec w}$. It follows that there exists $w'\in\mathcal{L}(u,v)$, $w'\prec w$, and a variable $x_j$ such that $x_1w=x_jw'$. Obviously, we have $j\geq2$. But this equality shows that $\nu_1(w')>\nu_1(w)$, which is impossible since $w'\prec w$.
\end{proof}

\begin{Lemma}\label{desc} Let $I=(\mathcal{L}(u,v))$ be a completely lexsegment ideal which has linear quotients with respect to the ordering $\prec$ of the generators. If $w\in\mathcal{L}(u,v)$ and $s\in\set(w)$, then
	\[g(x_sw)=
\left\{\begin {array}{cc}
			x_s w/x_1, & \mbox{if}\ x_sw\geq_{lex} x_1v,\\
			&\\
			x_s w/x_{\max(w)}, & \mbox{if}\ x_sw<_{lex} x_1v.
	\end{array}\right. \]
\end{Lemma}
\begin{proof} Let $u=x_1^{a_1}\ldots x_n^{a_n},\ v=x_1^{b_1}\ldots x_n^{b_n},\ a_1>0,$ and $w=x_1^{\alpha_1}\ldots x_n^{\alpha_n}$. 

In the first place we consider $$x_sw\geq_{lex}x_1v.$$ Since, by Lemma \ref{1nset}, we have $s\geq 2$, the above inequality shows that $\nu_1(w)\geq1$. We have to show that $g(x_sw)=x_s w/x_1$, that is $x_s w/x_1=\min_{\prec}\{w'\in\mathcal{L}(u,v)\ |\ x_sw\in I_{\preceq w'}\}$. It is clear that $v\leq_{lex} x_sw/x_1<_{lex}w\leq_{lex}u$, hence $x_sw/x_1\in\mathcal{L}(u,v)$. Let $w'\in\mathcal{L}(u,v)$ such that $x_sw\in I_{\preceq w'}$. We have to show that $x_sw/x_1\preceq w'$. Let $w''\in\mathcal{L}(u,v),\ w''\preceq w'$ such that $x_sw=w''x_j$, for some variable $x_j$. Then
$w''=x_sw/x_j\succeq x_sw/x_1
$
by the definition of our ordering $\prec$. This implies that $w'\succeq x_sw/x_1$.

Now we have to consider the second inequality, 
\begin{eqnarray}x_sw<_{lex}x_1v.\label{star}
\end{eqnarray}

Since $s\in\set(w)$, we have $x_sw\in I_{\prec w}$, that is there exists $w'\in\mathcal{L}(u,v)$, $w'\prec w$, and a variable $x_j,\ j\neq s$, such that 
\begin{eqnarray}x_sw=x_jw'.\label{1}
\end{eqnarray}

If $j=1$, then $x_sw=x_1w'\geq_{lex}x_1v$, contradiction. Hence $j\geq2$. We also note that $x_j|w$ since $j\neq s$, thus $j\leq\max(w)$. The following inequalities hold: 
\begin{eqnarray}
x_s w/x_{\max(w)}\geq_{lex} x_s w/x_j= w'\geq_{lex} v. \label{4}
\end{eqnarray}

If $\nu_1(w)<a_1$, we obviously get $x_sw/x_{\max(w)}\leq_{lex}u$. Let $\nu_1(w)=a_1$. From the inequality (\ref{star}) we obtain  $a_1\leq b_1+1$.  

If $a_1=b_1$ then $u=x_1^{a_1}x_2^{d-a_1}$ and $v=x_1^{a_1}x_n^{d-a_1}$ by Theorem \ref{completelylex}. Since $w\leq_{lex} u$, by using Lemma \ref{lema1}, we have
$x_sw/x_{\max(w)}\leq_{lex}x_su/x_{\max(u)}=x_su/x_2\leq_{lex} u,
$ the last inequality being true by Lemma \ref{1nset}. Therefore, $x_s w/x_{\max(w)}\in \mathcal{L}(u,v).$

If $a_1=b_1+1$ then the condition (c) in Theorem \ref{completelylex} holds. Let $z$ be the largest monomial with respect to the lexicographical order such that $z<_{lex}v$. Since $x_sw/x_1<_{lex}v$ by hypothesis, we also have $x_sw/x_1\leq_{lex}z$. By Lemma \ref{lema1} we obtain 
$x_sw/(x_1x_{\max(x_sw/x_1)})\leq_{lex}z/x_{\max(z)}.
$
Next we apply the condition (c) from Theorem \ref{completelylex} and get the following inequalities:

\begin{eqnarray}
x_1\frac{x_sw}{x_1x_{\max\left(\frac{x_sw}{x_1}\right)}}\leq_{lex}x_1\frac{z}{x_{\max(z)}}\leq_{lex}u.\label{2}
\end{eqnarray}
From the equality (\ref{1}) we have $w'=x_sw/x_j$.  As $j\neq1,\ \nu_1(w')=\nu_1(w),$ and the  inequality $w'\prec w$   gives $w'>_{lex} w,$ that is $x_s w/x_j>_{lex}w,$ which implies that $x_s>_{lex} x_j$. This shows that  $s<j\leq\max(w)$. Now looking at the inequalities (\ref{2}), we have
\begin{eqnarray}x_sw/x_{\max(w)}\leq_{lex}u.\label{3}
\end{eqnarray}
From (\ref{3}) and (\ref{4}) we obtain  $x_sw/x_{\max(w)}\in\mathcal{L}(u,v)$. 

It remains to show that $x_s w/x_{\max(w)}=\min_{\prec}\{w'\in\mathcal{L}(u,v)\mid x_sw\in I_{\preceq w'}\}.$ Let $\tilde{w}=\min_{\prec}\{w'\in\mathcal{L}(u,v)\mid x_sw\in I_{\preceq w'}\}.$ We obviously have $\tilde{w}\preceq x_s w/x_{\max(w)}\prec w.$ By the choice of $\tilde{w}$ we have $$x_sw=x_t\tilde{w}$$ for some variable $x_t.$ 

If $t=s$ we get $w=\tilde{w}$ which is impossible since $\tilde{w}\prec w.$ Therefore, $t\neq s.$
Then $x_t| w,$  so $t\leq\max(w)$. It follows that $\tilde{w}=x_sw/x_t\leq_{lex}x_sw/x_{\max(w)}$. If $t=1$ we have $x_1\tilde{w}=x_sw <_{lex}x_1 v,$ which implies that $\tilde{w}<_{lex}v,$ contradiction. Therefore $t\neq 1$ and, moreover, $\tilde{w}\succeq x_sw/x_{\max(w)}$, the  inequality being true by the definition of the ordering $\prec$. 
 This yields $\tilde{w}=x_sw/x_{\max(w)}.$ Therefore we have proved that $x_sw/x_{\max(w)}=g(x_sw)$.
\end{proof}

After this preparation, we  prove the following 

\begin{Theorem} Let $u=x_1^{a_1}\ldots x_n^{a_n},\ v=x_1^{b_1}\ldots x_n^{b_n},\ u,v\in\mathcal{M}_d$, with $u\geq_{lex}v$,  and $I=(\mathcal{L}(u,v))$ be a completely lexsegment ideal which has linear resolution. Then the decomposition function $g:M(I)\rightarrow G(I)$ associated to the ordering  $\prec$ of the generators from $G(I)$ is regular.
\end{Theorem}

\begin{proof} Let $w\in\mathcal{L}(u,v)$ and $s\in\set(w)$. We have to show that $\set(g(x_sw))\subset\set(w)$.

Let $t\in\set(g(x_sw))$. In order to prove that $t\in\set(w)$, that is $x_tw\in I_{\prec w}$, we will consider the following two cases:

\underline{\textit{Case 1}}: Let $x_sw\geq_{lex}x_1v$. By Lemma \ref{desc}, $g(x_sw)=x_sw/x_1$. Since $t\in\set(g(x_sw))$, we have $$\frac{x_tx_sw}{x_1}\in I_{\prec\frac{x_sw}{x_1}},$$ so there exists $w'\prec x_sw/x_1$, $w'\in\mathcal{L}(u,v)$, and a variable $x_j$, such that
$x_t x_sw/x_1=x_jw',
$
that is
	\begin{eqnarray}x_tx_sw=x_1x_jw'.\label{6}
\end{eqnarray}
By Lemma \ref{1nset}, $s,t\neq1$ and, since $w'\prec\ x_sw/x_1$, we have $j\neq t$. Note also that $w'\prec w$ since $\nu_1(w')<\nu_1(w)$. If $j=s$ then $x_tw=x_1w'\in I_{\prec w}$ and $t\in\set(w)$.

Now let $j\neq s$. If $j=1$, we have 
$x_tx_sw=x_1^2w',
$
which implies that $\nu_1(w')=\nu_1(w)-2$. The following inequalities hold:
$v<_{lex}x_1w'/x_s<_{lex}w\leq_{lex}u,
$
the first one being true since $v\leq_{lex}w'$, so $\nu_1(v)\leq\nu_1(w')$. These inequalities show that $x_1w'/x_s\in\mathcal{L}(u,v)$. But we also have $x_1w'/x_s\prec w$, hence $x_1w'/x_s\in I_{\prec w}$.

To finish this case we only need  to treat the case $j\neq1,\ j\neq s$. We are going to show that at least one of the monomials $x_1w'/x_s$ or $x_jw'/x_s$ belongs to $I_{\prec w}.$ In any case this will lead to the conclusion that $x_t w\in I_{\prec w}$ by using (\ref{6}).

From the equality (\ref{6}), we have $x_j|w$, hence $j\leq\max(w)$, and $\nu_1(w')=\nu_1(w)-1$. Since $w'\prec\ x_sw/x_1$ and $\nu_1(w')=\nu_1(w)-1=\nu_1(x_sw/x_1)$, we get 
\begin{eqnarray}w'>_{lex}x_sw/x_1,\label{7} 
\end{eqnarray} which gives
	\[x_1w'/x_s>_{lex}v.
\]
If the inequality \begin{eqnarray}x_1w'/x_s\leq_{lex}u\label{8}\end{eqnarray}holds, then we get $x_1w'/x_s\in\mathcal{L}(u,v)$. We also note that $\nu_1(x_1w'/x_s)=\nu_1(w)$ and $x_1w'/x_s>_{lex}w$ (by (\ref{7})). Therefore $x_1w'/x_s\prec w$ and we may write
$x_tw=x_j(x_1w'/x_s)\in I_{\prec w}.
$
This implies that $t\in\set(w)$.

Now we look at the monomial $x_jw'/x_s$ for which we have $\nu_1(x_jw'/x_s)=\nu_1(w')<\nu_1(w)$, so
$x_jw'/x_s<_{lex}w\leq_{lex}u.
$
If the inequality \begin{eqnarray}x_jw'/x_s\geq_{lex}v\label{9}\end{eqnarray}holds, we obtain  $x_jw'/x_s\in\mathcal{L}(u,v)$. Obviously we have $x_jw'/x_s\prec w$. By using (\ref{6}), we may write 
$x_tw=x_1(x_jw'/x_s)\in I_{\prec w},
$
which shows that $t\in\set(w)$.

To finish the proof in the Case $1$ we need to consider the situation when both inequalities (\ref{8}) and (\ref{9}) fail. Hence, let
	\[x_1w'/x_s>_{lex}u\ \mbox{and}\ x_jw'/x_s<_{lex}v.
\]
We will show that this inequalities cannot hold simultaneously. Comparing the exponents of $x_1$ in the monomials involved in the above inequalities, we obtain  $\nu_1(w')=b_1\geq a_1-1$. Since, by hypothesis, $x_sw>_{lex} x_1v,$ we have $\nu_1(w)>b_1.$ On the other hand, $w\leq_{lex} u$ implies that $\nu_1(w)\leq a_1.$ So $b_1=a_1-1$ and $\mathcal{L}(u,v)$ satisfies the condition (c) in Theorem \ref{completelylex}. Let, as usually, $z$ be the largest monomial with respect to the lexicographical order such that $z<_{lex}v$.

Since $x_jw'/x_s<_{lex}v$, we have $x_jw'/x_s\leq_{lex}z$. By Lemma \ref{lema1} and using the condition $x_1z/x_{\max(z)}\leq_{lex}u$, we obtain:
$x_1x_jw'/(x_sx_{\max(x_jw'/x_s)})\leq_{lex}u.
$
But our  assumption was that
$u<_{lex}x_1w'/x_s.
$
Therefore, combining the last two inequalities, after cancellation, one obtains that
$x_j<_{lex}x_{\max(x_jw'/x_s)}=x_{\max(x_tw/x_1)}=x_{\max(x_tw)}.
$
This leads to the inequality $j>\max(x_tw)$ and, since $j\leq\max(w)$, we get $\max(w)>\max(x_tw)$, which is impossible.

\underline{\textit{Case 2}}: Let $x_sw<_{lex}x_1v$. Then $g(x_sw)=x_sw/x_{\max(w)}$. In particular we have $x_sw/x_{\max(w)}\prec w$. Indeed, since $s\in \set(w),$ we have  $x_sw\in I_{\prec w}$, that is there exists $w'\in \mathcal{L}(u,v), w'\prec w,$ such that $x_sw\in I_{\preceq w'}.$ By the definition of the decomposition function we have $g(x_sw)\preceq w'$ and next we get  $g(x_sw)\prec w.$ Since $\nu_1(x_sw/x_{\max(w)})=\nu_1(w)$, the above inequality implies that $x_sw/x_{\max(w)}>_{lex}w$, that is $x_s>_{lex}x_{\max(w)}$ which means that $s<\max(w)$.

As $t\in\set(g(x_sw))$, there exists $w'\prec\ x_sw/x_{\max(w)}$, $w'\in\mathcal{L}(u,v)$, and a variable $x_j$, such that
	\[x_tx_sw/x_{\max(w)}=x_jw',
\]
that is
	\begin{eqnarray}x_tx_sw=x_jx_{\max(w)}w'.\label{10}
\end{eqnarray}
As in the previous case, we would like to show that one of the monomials $x_{\max(w)}w'/x_s$ or $x_jw'/x_s$ belongs to $\mathcal{L}(u,v)$ and it is strictly less than $w$ with respect to $\prec$. In this way we  obtain  $x_tw\in I_{\prec w}$ and $t\in\set(w)$.

We begin our proof noticing that $s,t\neq1$, by Lemma \ref{1nset}. The equality $j=t$ is impossible since $w'\neq\ x_sw/x_{\max(w)}$. If $j=s$, then $x_tw=w'x_{\max(w)}\in I_{\preceq w'}$. But $w'\prec\ x_sw/x_{\max(w)}\prec w$, hence $x_tw\in I_{\prec w}$. 

Let $j\neq s,t$. From the equality (\ref{10}) we have $x_j|w$, so $j\leq\max(w)$. We firstly consider $j=1$. Then the equality (\ref{10}) becomes
	\begin{eqnarray}x_tx_sw=x_1x_{\max(w)}w'.\label{11}
\end{eqnarray}
Since $s< \max(w)$, we have
$x_{\max(w)}w'/x_s<_{lex}w'\leq_{lex}u.
$
If the inequality $x_{\max(w)}w'/x_s$ $\geq_{lex}v$ holds too, then $x_{\max(w)}w'/x_s\in\mathcal{L}(u,v)$ and, as $\nu_1(w')<\nu_1(w)$, it follows that $x_{\max(w)}w'/x_s\prec w$. From (\ref{11}), we have $x_tw=x_1(x_{\max(w)}w'/x_s)\in I_{\prec w}$, hence $t\in\set(w)$. 

From the inequality $x_sw<_{lex}x_1v$, we get
	\[x_sw<_{lex}x_1w',
\]
so
	\[x_1w'/x_s>_{lex}w.
\]
Let us assume that $x_1w'/x_s\leq_{lex}u$. Since $\nu_1(x_1w'/x_s)=\nu_1(w)$, by using the definition of the ordering $\prec$ we get  $x_1w'/x_s\in I_{\prec w}$. Then we may write $x_tw=x_{\max(w)}(x_1w'/x_s)\in I_{\prec w}$.

It remains to consider that 
$x_{\max(w)}w'/x_s<_{lex}v\ \mbox{and}\ x_1w'/x_s>_{lex}u.
$
Proceeding as in the case 1 we show that we reach a contradiction and this ends the proof for $j=1$.
We only need to notice that we have to consider $b_1\leq a_1-1.$ Indeed, we can not have $b_1=a_1$ since one may find in $\mathcal{L}(u,v)$ at least two monomials, namely $w$ and $w',$ with $\nu_1(w')<\nu_1(w).$ 

Finally, let $j\neq1$. Recall that in the equality (\ref{10})  we have $j\neq 1,t,s $ and $s<\max(w)$. From (\ref{10}) we obtain  $\nu_1(w)=\nu_1(w')$.  Since $w'\prec \ x_s w/x_{\max(w)}$, we have $w'>_{lex} x_sw/x_{\max(w)}$, that is \begin{eqnarray}w'x_{\max(w)}>_{lex}x_sw.\label{12}\end{eqnarray}

Replacing $w'x_{\max(w)}$ by $x_tx_sw/x_j$ in (\ref{12}), we get $x_t>_{lex}x_j$, which means $t<j$. It follows that:
$x_{\max(w)}w'/x_s=x_tw/x_j>_{lex} w\geq_{lex}v.
$
Since $s<\max(w)$, as in the proof for $j=1$, we have $x_{\max(w)}w'/x_s\leq_{lex}u$. Therefore $x_{\max(w)}w'/x_s\in\mathcal{L}(u,v)$. In addition, from (\ref{12}), $x_{\max(w)}w'/x_s>_{lex}w$ and $\nu_1(x_{\max(w)}w'/x_s)=\nu_1(w)$, so $x_{\max(w)}w'/x_s\prec w$. In other words, we have got that
$x_tw=x_j(x_{\max(w)}w'/x_s)\in I_{\prec w}
$
and $t\in\set(w)$.
\end{proof}

The general problem of determining the resolution of arbitrary lexsegment ideals is not completely solved. The resolutions of the lexsegment ideals with linear quotients are described in \cite{HT} using iterated mapping cones. We recall this construction from \cite{HT}. Suppose that the monomial ideal $I$ has linear quotients with respect to the ordering $u_1,\ldots,u_m$ of its minimal generators. Set $I_j=(u_1,\ldots,u_j)$ and $L_j=(u_1,\ldots,u_j):u_{j+1}.$ Since $I_{j+1}/I_j\simeq S/L_j,$ we get the exact sequences $$0\rightarrow S/L_j \rightarrow S/I_j\rightarrow S/I_{j+1}\rightarrow 0,$$ where the morphism $S/L_j \rightarrow S/I_j$ is the multiplication by $u_{j+1}.$ Let $F^{(j)}$ be a graded free resolution of $S/I_j$, $K^{(j)}$ the Koszul complex associated to the regular sequence $x_{k_1},\ldots,x_{k_l}$ with $k_i\in\set(u_{j+1}),$ and $\psi^{(j)}:K^{(j)}\rightarrow F^{(j)}$ a graded complex morphism lifting the map $S/L_j \rightarrow S/I_j$. Then the mapping cone $C(\psi^{(j)})$ of $\psi^{(j)}$ yields a free resolution of $S/I_{j+1}.$ By iterated mapping cones we obtain step by step a graded free resolution of $S/I.$

\begin{Lemma}\cite{HT} Suppose $\deg\ u_1 \leq \deg\ u_2 \leq \ldots\leq \deg\ u_m.$ Then the iterated mapping
cone $\mathbb{F}$, derived from the sequence $u_1,\ldots,u_m,$ is a minimal graded free resolution
of $S/I$, and for all $i > 0$ the symbols
\[f(\sigma; u)\ \mbox{with}\ u\in G(I),\ \sigma \subset \set(u),\  |\sigma| = i - 1
\]
form a homogeneous basis of the $S-$module $F_i$. Moreover $\deg(f(\sigma; u)) = |\sigma| +\deg(u)$.
\end{Lemma} 

\begin{Theorem}\cite{HT} Let $I$ be a monomial ideal of $S$ with linear quotients, and $\mathbb{F}_{\bullet}$ the graded
minimal free resolution of $S/I$. Suppose that the decomposition function $g : M(I) \rightarrow G(I)$ is regular. Then the chain map $\partial$ of $\mathbb{F}_{\bullet}$ is given by
\[\partial(f(\sigma; u)) = -\sum_{s\in\sigma}(-1)^{\alpha(\sigma;s)}x_sf(\sigma\setminus s;u)+\sum_{s\in\sigma}(-1)^{\alpha(\sigma;s)}\frac{x_su}{g(x_su)}f(\sigma\setminus s;g(x_su)),\]
if $\sigma\neq\emptyset$, and
\[\partial(f(\emptyset; u)) = u\] otherwise.
Here $\alpha(\sigma;s)=|\{t\in\sigma\ |\ t<s\}|$.
\end{Theorem}

In our specific context we get the following
\begin{Corollary} Let $I=(\mathcal{L}(u,v))\subset S$ be a completely lexsegment ideal with linear quotients with respect to $\prec$ and $\mathbb{F}_{\bullet}$ the graded minimal free resolution of $S/I$. Then the chain map of $\mathbb{F}_{\bullet}$ is given by
\[\partial(f(\sigma; w)) = -\sum_{s\in\sigma}(-1)^{\alpha(\sigma;s)}x_sf(\sigma\setminus s;w)+\sum_{\stackrel{s\in\sigma:}{x_sw\geq_{lex}}x_1v}(-1)^{\alpha(\sigma;s)}x_1f\left(\sigma\setminus s;\frac{x_sw}{x_1}\right)+\]\[+\sum_{\stackrel{s\in\sigma:}{x_sw<_{lex}}x_1v}(-1)^{\alpha(\sigma;s)}x_{\max(w)}f\left(\sigma\setminus s;\frac{x_sw}{x_{\max(w)}}\right),\]
if $\sigma\neq\emptyset$, and
\[\partial(f(\emptyset; w)) = w\] otherwise. For convenience we set $f(\sigma;w)=0$ if $\sigma\nsubseteq\set{w}$.
\end{Corollary}

\begin{Example}\rm
Let $u=x_1^2x_2$ and $v=x_2^3$  be monomials in the polynomial ring $S=k[x_1,x_2,x_3]$. Then $$\mathcal{L}(u,v)=\{x_2^3,\ x_1x_2^2,\ x_1x_2x_3,\ x_1x_3^2,\ x_1^2x_2\}.$$ The ideal $I=(\mathcal{L}(u,v))$ is a completely lexsegment ideal with linear quotients with respect to this ordering of the generators. We denote $u_1=x_2^3,\ u_2=x_1x_2^2,\ u_3=x_1x_2x_3,\ u_4=x_1x_3^2,\ u_5=x_1^2x_2$. We have  $\set(u_1)=\emptyset,\ \set(u_2)=\{2\},\ \set(u_3)=\{2\},\ \set(u_4)=\{2\},\ \set(u_5)=\{2,3\}$. Let $\mathbb{F}_{\bullet}$ be the minimal graded free resolution of $S/I$.

Since $\max\{|\set(w)|\mid w\in\mathcal{L}(u,v)\}=2$, we have $F_i=0$, for all $i\geq4$.

A basis for the $S-$module $F_1$ is $\{f(\emptyset;u_1),\ f(\emptyset;u_2),\ f(\emptyset;u_3),\ f(\emptyset;u_4),\ f(\emptyset;u_5)\}$.
	
A basis for the $S-$module $F_2$ is $$\{f(\{2\};u_2),\ f(\{2\};u_3),\ f(\{2\};u_4),\ f(\{2\};u_5),\ f(\{3\};u_5)\}.$$

A basis for the $S-$module $F_3$ is $\{f(\{2,3\};u_5)\}$.

We have the minimal graded free resolution $\mathbb{F}_{\bullet}$:
	\[0\rightarrow S(-5)\stackrel{\partial_2}{\rightarrow} S(-4)^5\stackrel{\partial_1}{\rightarrow} S(-3)^5\stackrel{\partial_0}{\rightarrow} S\rightarrow S/I\rightarrow 0
\]
where the maps  are
	\[\partial_0(f(\emptyset;u_i))=u_i,\ \mbox{for}\ 1\leq i\leq 5, 
\]
so \[\partial_0=
\left(\begin {array}{ccccc}
			x_2^3& x_1x_2^2& x_1x_2x_3& x_1x_3^2& x_1^2x_2
	\end{array}\right) .\]
\[\begin{array}{lll}
\partial_1(f(\{2\};u_2))&= &-x_2f(\emptyset;u_2)+x_1f(\emptyset;u_1),\\
\partial_1(f(\{2\};u_3))&= &-x_2f(\emptyset;u_3)+x_3f(\emptyset;u_2),\\
\partial_1(f(\{2\};u_4))&= &-x_2f(\emptyset;u_4)+x_3f(\emptyset;u_3),\\
\partial_1(f(\{2\};u_5))&= &-x_2f(\emptyset;u_5)+x_1f(\emptyset;u_2),\\
\partial_1(f(\{3\};u_5))&= &x_3f(\emptyset;u_5)-x_1f(\emptyset;u_3),\\
\end{array}\]
so
\[\partial_1=
\left(\begin {array}{ccccc}
			x_1& 0& 0& 0& 0\\
			-x_2& x_3& 0& x_1& 0\\
			0& -x_2& x_3& 0& -x_1\\
			0& 0& -x_2& 0& 0\\
			0& 0& 0& -x_2& x_3
	\end{array}\right) .\]

	\[\partial_2(f(\{2,3\};u_5))=-x_2f(\{3\};u_5)+x_3f(\{2\};u_5)+x_1f(\{3\};u_2)-x_1f(\{2\};u_3)=\]\[=-x_2f(\{3\};u_5)+x_3f(\{2\};u_5)-x_1f(\{2\};u_3),
\]
since $\{3\}\nsubseteq\set(u_2)$, so\[\partial_2=
\left(\begin {array}{c}
			0\\
			-x_1\\
			0\\
			x_3\\
			-x_2
	\end{array}\right) .\]
\end{Example}
\section{Non-completely lexsegment ideals with linear resolutions}

\begin{Theorem}\label{noncompletely}
Let $u=x_1^{a_1}\ldots x_n^{a_n},\ v=x_2^{b_2}\ldots x_n^{b_n}$ be monomials of degree $d$ in $S,$  $a_1> 0.$ Suppose that the ideal $I=({\mathcal L}(u,v))$ is not completely lexsegment ideal. Then $I$ has linear resolution if and only if $I$ has linear quotients.
\end{Theorem}

\begin{proof}
We only have to prove that if $I$ has linear resolution then $I$ has linear quotients for a suitable ordering of its minimal monomial
 generators. By \cite[Theorem 2.4]{ADH}, since $I$ has linear resolution,  $u$ and $v$ have the form: $$u=x_1x_{l+1}^{a_{l+1}}\ldots
  x_n^{a_n},\ v=x_lx_n^{d-1}, \text {for some } l\geq 2.$$ Then the ideal $I=({\mathcal L}(u,v))$ can be written as a sum of ideals $I=J+K,$ where $J$ is the ideal generated by all the monomials of ${\mathcal L}(u,v)$ which are not divisible by $x_1$ and $K$ is
 generated by all the monomials of ${\mathcal L}(u,v)$ which are  divisible by $x_1.$  More precise, we have $$J=(\{w\mid
  x_2^d\geq_{lex} w\geq_{lex} v\})$$ and $$K=(\{w\mid u\geq_{lex} w\geq_{lex} x_1x_n^{d-1}\}).$$ One may see that $J$ is generated by
 the initial lexsegment $\mathcal{L}^i(v)\subset k[x_2,\ldots,x_n],$ and hence it has linear quotients with respect to
 lexicographical order $>_{lex}.$ Let $G(J)=\{g_1\prec \ldots\prec g_m\},$ where $g_i\prec g_j\ \text{if and only if}\
 g_i>_{lex} g_j.$ The ideal $K$ is isomorphic with the ideal generated by the final lexsegment of degree $d-1$
 $$\mathcal{L}^f(u/x_1)=\{w\mid u/x_1\geq_{lex} w\geq_{lex} x_n^{d-1},\ \deg(w)=d-1\}.$$ Since final lexsegments
  are stable with respect to the order  $x_n>\ldots >x_1$ of the variables, it follows that the ideal $K$ has linear quotients
 with respect to $>_{\overline{lex}},$ where by $\overline{lex}$ we mean the lexicographical order corresponding to 
 $x_n>\ldots >x_1.$ Let $G(K)=\{h_1\prec \ldots\prec h_p\},$ where $h_i\prec h_j\ \text{if and only if}\
 h_i>_{\overline{lex}} h_j.$ We consider the following ordering of the monomials of $G(I):$ $$G(I)=\{g_1\prec \ldots\prec
 g_m\prec h_1\prec \ldots\prec h_p\}.$$ We claim that, for this ordering of its minimal monomial generators, $I$ has linear quotients. In order to check this, we firstly notice that $I_{\prec g}:g=J_{\prec g}:g$ for every  $g\in G(J).$ Since $J$
 has linear quotients with respect to $\prec$ it follows that $J_{\prec g}:g$ is generated by variables. Now it is enough
to show that, for any generator $h$ of $K,$ the colon ideal $I_{\prec h}:h$ is generated by variables. 
We note that $$I_{\prec h}:h=J:h + K_{\prec h}:h.$$ Since $K$ is with linear quotients, we already know that
$K_{\prec h}:h$ is generated by variables. Therefore we only need to prove that $J:h$ is generated by variables. We
will show that $J:h=(x_2,\ldots,x_l)$ and this will end our proof. Let $m\in J:h$ be a monomial. It follows that $mh\in J.$ Since $h$ is a generator of $ K,$ $h$ is of the form $h=x_1x_{l+1}^{\alpha_{l+1}}\ldots x_n^{\alpha_n},$ that is $h\not\in (x_2,\ldots,x_l).$ But this implies that $m$ must be in the ideal $(x_2,\ldots,x_l).$ For the reverse inclusion, let $2\leq t\leq  l.$ Then $x_t h=x_1 \gamma$ for some monomial $\gamma,$ of degree $d.$ Replacing $h$ in the equality we get  $\gamma=x_tx_{l+1}^{\alpha_{l+1}}\ldots x_n^{\alpha_n}$ which shows that $\gamma$ is a generator of $J.$ Hence $x_th\in J.$
\end{proof}

\begin{Example}\rm 
Let $I=(\mathcal{L}(u,v))\subset k[x_1,\ldots,x_6]$ be the lexsegment ideal of degree $4$ determined by the monomials $u=x_1x_3^2x_5$ and $v=x_2x_6^3.$  $I$ is not a completely lexsegment ideal as it follows applying \cite[Theorem 2.3]{DH}, but $I$ has linear resolution by \cite[Theorem 2.4]{ADH}. $I$ has linear quotients if we order its minimal monomial generators as indicated in the proof of the above theorem. On the other hand, if we order the generators of $I$ using the order relation defined in the proof of Theorem  \ref{colex} we can easy see that $I$ does not have linear quotients. Indeed, following the definition of the order relation from Theorem \ref{colex} we should take $$G(I)=\{x_2^4\prec x_2^3x_3\prec\ldots\prec x_2x_6^3\prec x_1x_3^2x_5\prec x_1x_3^2x_6\prec x_1x_3x_4^2\prec\ldots\prec x_1x_6^3\}.$$ For $h=x_1x_3x_4^2$ one may easy check that  $I_{\prec h}:h$ is not generated by variables.
\end{Example}

\begin{Example}\rm Let $u=x_1x_3^2$, $v=x_2x_4^2$ be monomials in $k[x_1,\ldots, x_4]$. Then $I=(L(u,v))\subset k[x_1,\ldots, x_4]$ is a non-completely ideal with linear resolution and, by the proof of Theorem \ref{noncompletely}, $I$ has linear quotients with respect to the following ordering  of its minimal  monomial generators:
$$x_2^3,\ x_2^2x_3,\ x_2^2x_4,\ x_2x_3^2,\ x_2x_3x_4,\ x_2x_4^2,\ x_1x_4^2,\ x_1x_3x_4,\ x_1x_3^2.$$
We note that $\set(x_1x_4^2)=\{2\}$ and $\set(g(x_1x_2x_4^2))=\set(x_2x_4^2)=\{2,3\}\nsubseteq\set(x_1x_4^2)$, so the decomposition function is not regular for this ordering of the generators.
\end{Example}
\section{Cohen-Macaulay lexsegment ideals} \label{Section3}

In this section we study the dimension and the depth of arbitrary lexsegment ideals. These results are applied to describe the lexsegments ideals which are Cohen-Macaulay. We begin with the study of the dimension. As in the previous sections, let $d\geq 2$ be an integer. We denote $\frak{m}=(x_1,\ldots,x_n).$ It is clear that if $I=(\mathcal{L}(u,v))\subset S$ is a lexsegment ideal of degree $d$ then $\dim(S/I)=0$ if and only if $I=\frak{m}^d.$

\begin{Proposition}\label{dim} Let $u=x_1^{a_1}\ldots x_n^{a_n},\ v=x_q^{b_q}\ldots x_n^{b_n}$, $1\leq q\leq n$, $a_1,b_q>0,$ be two monomials of degree $d$ such that $u\geq_{lex}v$ and let $I$ be the lexsegment ideal generated by $\mathcal{L}(u,v)$. We assume that $I\neq\frak{m}^d$. Then
\[\dim(S/I)=
\left\{\begin {array}{ll}
			n-q, & \mbox{if}\ 1\leq q<n,\\
			1, & \mbox{if}\ q=n.
	\end{array}\right. \]
\end{Proposition}

\begin{proof} For $q=1$, we have $I\subset(x_1)$. Obviously $(x_1)$ is a minimal prime of $I$ and $\dim(S/I)=n-1$.

Let $q=n$, that is $v=x_n^d$ and $\mathcal{L}(u,v)=\mathcal{L}^f(u)$. We may write the ideal $I$ as a sum of two ideals, 
$I=J+K,$ where
$J=(x_1\mathcal{L} (u/x_1,x_n^{d-1}))
$
and
$K=(\mathcal{L}(x_2^d,x_n^d)).
$
Let $p\supset I$ be a monomial prime ideal. If $x_1\in p$, then $J\subseteq p$. Since $p$ also contains $K$, we have  $p\supset (x_2,\ldots,x_n)$. Hence $p=(x_1,x_2,\ldots,x_n)$. If $x_1\notin p$, we obtain  $(x_2,\ldots,x_n)\subset p$. Hence, the only minimal prime ideal of $I$ is $(x_2,\ldots,x_n)$. Therefore, $\dim(S/I)=1$.

Now we consider $1<q<n$ and write $I$ as before,
$I=J+K,
$
where $J=(x_1\mathcal{L}(u/x_1,x_n^{d-1}))$ and $K=(\mathcal{L}(x_2^d,v))$. 

Firstly we consider $u=x_1^d.$ Let $p\supset I$ be a monomial prime ideal. Then $p\ni x_1$ and, since $p\supset K,$ we also have $p\supset (x_2,\ldots,x_q).$ Hence $(x_1,\ldots,x_q)\subset p.$ Since $I\subset (x_1,\ldots,x_q),$ it follows that $(x_1,\ldots,x_q)$ is the only minimal prime ideal of $I.$ Therefore $\dim(S/I)=n-q.$

Secondly, let $a_1>1$ and $u\neq x_1^d.$ The lexsegment $\mathcal{L}(u/x_1,x_n^{d-1})$ contains the lexsegment $\mathcal{L}(x_2^{d-1},x_n^{d-1})$. Let $p$ be a monomial prime ideal which contains $I$ and such that $x_1\not\in p$. Then $p\supset\mathcal{L}(x_2^{d-1},x_n^{d-1})$ which implies that $(x_2,\ldots,x_n)\subset p$. Obviously we also have $I\subset(x_2,\ldots,x_n)$, hence $(x_2,\ldots,x_n)$ is a minimal prime ideal of $I$.

Let $p\supset I$ be a monomial prime ideal which contains $x_1$. Since $p\supset K$, we also have $(x_2,\ldots,x_q)\subset p$. This shows that $(x_1,\ldots, x_q)$ is a minimal prime ideal of $I$. In conclusion, for $a_1>1$, the minimal prime ideals of $I$ are $(x_1,\ldots,x_q)$ and $(x_2,\ldots, x_n)$. Since $q\leq n-1$, we get  $\height(I)=q$ and $\dim(S/I)=n-q$.

Finally, let $a_1=1$, that is $u=x_1 x_l^{a_l}\ldots x_n^{a_n}$, for some $a_l>0$, $l\geq2$. As in the previous case, we obtain  $(x_1,\ldots,x_q)$  a minimal prime ideal of $I$. Now we look for those minimal prime ideals of $I$ which do not contain $x_1$.

If $a_l=d-1$, the ideal $J=(x_1\mathcal{L}(u/x_1,x_n^{d-1}))$ becomes $J=(x_1\mathcal{L}(x_l^{d-1},x_n^{d-1}))$. If $p\supset I$ is a monomial prime ideal such that  $x_1\notin p$, we get  $(x_l,\ldots,x_n)\subset p$, and, since $p$ contains $K$, we obtain  $(x_2,\ldots,x_q)\subset p$. This shows that if $q<l$ then $(x_2,\ldots,x_q,x_l,\ldots,x_n)$ is a minimal prime ideal of $I$ of height $q+n-l\geq q$, and if $q\geq l$, then $(x_2,\ldots,x_n)$ is a minimal prime ideal of height $n-1\geq q$. In both cases we may draw the conclusion that $\height(I)=q$ and, consequently, $\dim(S/I)=n-q$.

The last case we have to consider is $a_l<d-1$. Then $l<n$ and, with similar arguments as above, we obtain  $\dim(S/I)=n-q$.
\end{proof}

In order to study the depth of arbitrary lexsegment ideals, we  note that one can restrict to those lexsegments defined by monomials of the form $u=x_1^{a_1}\ldots x_n^{a_n},\ v=x_1^{b_1}\ldots x_n^{b_n}$ of degree $d$ with $a_1>0$ and $b_1=0$.

Indeed, if $a_1=b_1$, then $I=(\mathcal{L}(u,v))$ is isomorphic, as an $S-$module, with the ideal generated by the lexsegment $\mathcal{L}(u/x_1^{a_1},v/x_1^{b_1})$ of degree $d-a_1$. This lexsegment may be studied in the polynomial ring in a smaller number of variables.

If $a_1>b_1$, then $I=(\mathcal{L}(u,v))$ is isomorphic, as an $S-$module, with the ideal generated by the lexsegment $\mathcal{L}(u',v')$, where $u'=u/x_1^{b_1}$ has $\nu_1(u')=a_1-b_1>0$ and $v'=v/x_1^{b_1}$ has $\nu_1(v')=0$.

Taking into account these remarks, from now on, we  consider lexsegment ideals of ends $u=x_1^{a_1}\ldots x_n^{a_n}$, $v=x_q^{b_q}\ldots x_n^{b_n}$, for some $q\geq2$, $a_1,b_q>0$.

The first step in the depth's study  is the next

\begin{Proposition}\label{depthzero} Let $I=(\mathcal{L}(u,v))$, where $u=x_1^{a_1}\ldots x_n^{a_n}$, $v=x_q^{b_q}\ldots x_n^{b_n}$, $q\geq2$, $a_1,b_q>0$. Then $\depth(S/I)=0$ if and only if $x_nu/x_1\geq_{lex} v$.
\end{Proposition}
\begin{proof} Let $x_nu/x_1\geq_{lex} v$. We claim that
$(I\colon (u/x_1))=(x_1,\ldots,x_n).
$
Indeed, for $1\leq j\leq n$, the  inequalities 
$u\geq_{lex}x_ju/x_1\geq_{lex}x_nu/x_1\geq_{lex}v
$ hold.
They show that $x_ju/x_1\in I$ for $1\leq j\leq n$. Therefore $(x_1,\ldots,x_n)\subseteq(I\colon (u/x_1))$. The other inclusion is obvious. We conclude that $(x_1,\ldots,x_n)\in \Ass(S/I)$, hence $\depth(S/I)=0$.

For the converse, let us assume, by contradiction, that $x_nu/x_1<_{lex}v$. We will show that $x_1-x_n$ is regular on $S/I$. This will imply that $\depth(S/I)>0$, which contradicts our hypothesis. We firstly notice that, from the above inequality, we have $a_1-1=0$, that is $a_1=1$. Therefore, $u$ is of the form $u=x_1x_l^{a_l}\ldots x_n^{a_n}$, $l\geq2,\ a_l>0$. Moreover, we have $l\geq q.$

Let us suppose that $x_1-x_n$ is not regular on $S/I$, that is there exists at least a polynomial $f\notin I$ such that $f(x_1-x_n)\in I$. One may assume that all monomials of $\supp(f)$ do not belong to $I$. Let us choose such a polynomial $f=c_1w_1+\ldots+c_tw_t,\ c_i\in k,\ 1\leq i\leq t$, with $w_1>_{lex}w_2>_{lex}\ldots>_{lex}w_t,\ w_i\notin I,\ 1\leq i\leq t$.

Then $\ini_{lex}((x_1-x_n)f)=x_1w_1\in I$. It follows that there exists $\alpha\in G(I)$ such that
\begin{eqnarray}x_1w_1=\alpha\cdot\alpha\,'.\label{3.1}
\end{eqnarray}
for some monomial $\alpha\,'$. We have $x_1\nmid\alpha\,'$ since, otherwise, $w_1\in I$, which is false. Hence $\alpha$ is a minimal generator of $I$ which is divisible by $x_1$, that is $\alpha$ is of the form $\alpha=x_1\gamma$, for some monomial $\gamma$ such that $x_n^{d-1}\leq_{lex}\gamma\leq_{lex}u/x_1$. Looking at $(\ref{3.1})$, we get  $w_1=\gamma \alpha\,'$. This equality shows that $x_1\nmid w_1$. We claim that the monomial $x_nw_1$ does not cancel in the expansion of $f(x_1-x_n)$. Indeed, it is clear that $x_nw_1$ cannot cancel by some monomial $x_nw_i$,  $i\geq2$. But it also cannot cancel by some monomial of the form $x_1w_i$ since $x_nw_1$ is not divisible by $x_1$. Now we may draw the conclusion that there exists a monomial $w\notin I$ such that $w(x_1-x_n)\in I$, that is $wx_1,\ wx_n\in I$.

Let $w\notin I$ be a monomial  such that $wx_1,\ wx_n\in I$, let $\alpha,\ \beta\in\mathcal{L}(u,v)$ and $\alpha\,',\ \beta\,'$ monomials such that
\begin{eqnarray}x_1w=\alpha\cdot\alpha\,'\label{3.2}
\end{eqnarray}
and
\begin{eqnarray}x_nw=\beta\cdot\beta\,'.\label{3.3}
\end{eqnarray}
As before, we get  $x_1\nmid w$, hence $\beta$ must be a minimal generator of $I$ such that
$x_2^d\geq_{lex}\beta\geq_{lex}v.
$ By using (\ref{3.3}), we can see that $x_n$ does not divide $\beta\,',$ hence $x_n|\beta.$  It follows that $w$ is divisible by $\beta/ x_n$. $w$ is also divisible by $\alpha/x_1$. Therefore,
 $\delta=\lcm(\alpha/x_1,\beta/x_n) | w.$ If $\deg \delta\geq d$ there exists a variable $x_j,$ with $j\geq 2,$ such that $(x_j\beta/x_n )| \delta,$ thus $(x_j\beta/x_n) | w.$ It is obvious that $x_2^d\geq_{lex} x_j\beta/x_n \geq_{lex} \beta\geq_{lex} v,$ hence $x_j\beta/x_n $ is a minimal generator of $I$ which divides $w,$ contradiction. This implies that $\delta$ has the degree $d-1.$ This yields $\alpha/x_1=\beta/x_n.$ Then $\beta=x_n\alpha /x_1\leq_{lex} x_n u/x_1<_{lex} v,$ contradiction.

\end{proof}

\begin{Corollary}
Let $I=(\mathcal{L}(u,v))$, where $u=x_1^{a_1}\ldots x_n^{a_n}$, $v=x_q^{b_q}\ldots x_n^{b_n}$, $q\geq2$, $a_1,b_q>0$. Then $\projdim(S/I)=n$ if and only if $x_nu/x_1\geq_{lex} v$.
\end{Corollary}

\begin{Corollary}Let $I=\left(\mathcal{L}^f(u)\right)$ be the ideal generated by the final lexsegment defined by $u=x_1^{a_1}\ldots x_n^{a_n},\ a_1>0$. Then $\depth(S/I)=0$.
\end{Corollary}

\begin{Corollary}
Let $I=\left(\mathcal{L}^i(v)\right)$ be the ideal generated by the initial lexsegment defined by the monomial $v$. Then $\depth(S/I)=0$ if and only if $v\leq_{lex}x_1^{d-1}x_n$.
\end{Corollary}

Next we are going to characterize the lexsegment ideals $I$ such that $\depth\ S/I>0$, that is $x_nu/x_1<_{lex}v$, which implies that $u$ has the form $u=x_1x_l^{a_l}\ldots x_n^{a_n}$, for some $l\geq2,\ a_l>0$ and $l>q,$ or $l=q$ and $a_q\leq b_q$. We  denote $u\,'=u/x_1=x_l^{a_l}\ldots x_n^{a_n}$. Then we have $x_nu\,'<_{lex}v$. From the proof of Proposition \ref{depthzero} we know that $x_1-x_n$ is regular on $S/I$. Therefore 
	\[\depth(S/I)=\depth(S\,'/I\,')+1,
\]
where $S\,'=k[x_2,\ldots,x_n]$ and $I\,'$ is the ideal of $S\,'$ whose minimal monomial generating set is $G(I\,')=x_n\mathcal{L}(u\,',x_n^{d-1})\cup\mathcal{L}^i(v)$.

\begin{Lemma}\label{depthneq0} In the above notations and  hypothesis on the lexsegment ideal $I$, the following statements hold:
\begin{itemize}
	\item[(a)] If $v=x_2^d$ and $l\geq4$, then $\depth(S\,'/I\,')=l-3$.
	\item[(b)] If $v=x_2^{d-1}x_j$ for some $3\leq j\leq n-2$ and $l\geq j+2$ then $\depth(S\,'/I\,')=l-j-1$.
	\item[(c)] $\depth(S\,'/I\,')=0$ in all the other cases. 
\end{itemize}
\end{Lemma}

\begin{proof} (a) Let $v=x_2^d$ and $l\geq4$. The ideal $I\,'\subset S '$ is minimally generated by all the monomials $x_n\gamma$, where $x_n^{d-1}\leq_{lex}\gamma\leq_{lex}u\,'$, $\deg(\gamma)=d-1$, and by the monomial $x_2^d$. Then it is clear that $\{x_3,\ldots,x_{l-1}\}$ is a regular sequence on $ S\,'/I\,'$, hence
$$\depth S\,'/I\,'=\depth\frac{S\,'/I\,'}{(x_3,\ldots,x_{l-1})S\,'/I\,'}+l-3.
$$We have
$$\frac{S\,'/I\,'}{(x_3,\ldots,x_{l-1})S\,'/I\,'}\cong\frac{k[x_2,x_l,\ldots,x_n]}{I\,'\cap k[x_2,x_l,\ldots,x_n]}.
$$In this way we may reduce the computation of $\depth(S\,'/I\,')$ to the case (c).

(b) Let $v=x_2^{d-1}x_j$, for some $3\leq j\leq n-2$ and $l\geq j+2$. Hence $I\,'$ is minimally generated by the following set of monomials
	\[\{x_n\gamma\ |\ \gamma\ \mbox{monomial of degree}\ d-1\ \mbox{such that}\ x_n^{d-1}\leq_{lex}\gamma\leq_{lex}u\,' \}\cup\]\[\cup\{x_2^d,\ x_2^{d-1}x_3,\ldots,x_2^{d-1}x_j\}.
\]
Then $\{x_{j+1},\ldots,x_{l-1}\}$ is a regular sequence on $S\,'/I\,'$ and 
$$\depth S\,'/I\,'=\depth\frac{S\,'/I\,'}{(x_{j+1},\ldots,x_{l-1})S\,'/I\,'}+(l-j-1).
$$
Since
$$\frac{S\,'/I\,'}{(x_{j+1},\ldots,x_{l-1})S\,'/I\,'}\cong\frac{k[x_2,\ldots,x_j,x_l,\ldots,x_n]}{I\,'\cap k[x_2,\ldots,x_j,x_l,\ldots,x_n]},
$$
we may reduce the computation of $\depth(S\,'/I\,')$ to the case (c).

(c) In each of the cases that it remains to treat, we will show that $(x_2,\ldots,x_n)\in \Ass(S '/I ')$, that is there exists a monomial $w\notin I\,'$ such that $I\,':w=(x_2,\ldots,x_n)$. This implies that $\depth(S\,'/I\,')=0$.

{\textit{Subcase $C_1$}:} $v=x_2^d,\ l=2$. Then $w=x_n^{d-1}\notin I\,'$ and 
$x_n^{d-1}\leq_{lex}x_jw/x_n=x_jx_n^{d-2}\leq_{lex}x_2x_n^{d-2}\leq_{lex}x_l^{a_l}\ldots x_n^{a_n}=u',
$
for all $2\leq j\leq n$. Hence $\gamma=x_jw/x_n$ has the property that $x_n\gamma\in G(I\,')$. Therefore, $x_j\in I\,':w$ for all $2\leq j\leq n$. It follows that $I\,':w=(x_2,\ldots,x_n)$.

{\textit{Subcase $C_2$}:} $v=x_2^d,\ l=3$. Then $w=x_2^{d-1}x_n^{d-1}\notin I\,'$. Indeed, $x_2^d\nmid w$ and if we assume that there exists $x_n^{d-1}\leq_{lex}\gamma\leq_{lex}u\,'$, $\deg\gamma=d-1,$ such that $x_n\gamma|w$,  we obtain  $x_n\gamma|x_n^{d-1}$ which is impossible.

We show that $x_jw\in I\,'$ for all $2\leq j\leq n$. Indeed, $x_2w=x_2^dx_n^{d-1}\in I\,'$. Let $3\leq j\leq n$. Then $x_n^{d-1}\leq_{lex}x_jx_n^{d-2}\leq_{lex}x_3x_n^{d-2}\leq_{lex} u\,'$. It follows that $\gamma=x_jx_n^{d-2}$ has the property that $x_n\gamma=x_jx_n^{d-1}\in G(I\,')$. Since $x_n\gamma|x_jw$, we have  $x_jw\in I\,'$. This arguments shows that $I\,':w=(x_2,\ldots,x_n)$.

{\textit{Subcase $C_3$}:} $v=x_2^{d-1}x_j$ for some $3\leq j\leq n-1$ and $2\leq l\leq j+1$. Let us consider again the monomial $w=x_2^{d-1}x_n^{d-1}$. It is clear that $x_tw\in I$ for all $2\leq t\leq j$. Let $t\geq j+1$. Then $x_tw$ is divisible by $x_tx_n^{d-1}$. Since $x_tx_n^{d-2}$ satisfies the  inequalities
$x_n^{d-1}\leq_{lex}x_tx_n^{d-2}\leq_{lex}u\,',
$
we have  $x_tx_n^{d-1}\in G(I\,')$. It follows that $x_t w\in I\,'$ for $t\geq j+1$. Assume that $w\in I\,'$. Since $x_2^{d-1}x_t\nmid w$ for $2\leq t\leq j$, we should have $x_n\gamma|w$ for some $\gamma$ of degree $d-1$ such that $x_n^{d-1}\leq_{lex}\gamma\leq_{lex}u\,'$. Since $\gamma|x_2^{d-1}x_n^{d-2}$ and $\gamma\leq_{lex}u\,'$, we get  $l=2$ and $a_2=\nu_2(u\,')\geq\nu_2(\gamma)$. Let $\gamma=x_2^ax_n^{d-1-a}$, for some $a\geq1$. In this case we change the monomial $w.$ Namely, we consider the monomial $w'=x_2x_n^{d-2}$ which does not belong to $G(I\,')$ since it has degree $d-1$.

 If $a_2\geq2$, for any $j$ such that $2\leq j\leq n$, we have 
$x_n^{d-1}<_{lex}x_jw'/x_n=x_2x_jx_n^{d-3}<_{lex}x_l^{a_l}\ldots x_n^{a_n}=u\,'.
$
This shows that $x_jw'\in I\,'$ for $2\leq j\leq n$ and hence, $I\,':w=(x_2,\ldots,x_n)$.

If $a_2=1$, we take $w\,''=x_n^{d-1}\notin I\,'$. For all $j$ such that $2\leq j\leq n$, we have 
$x_n^{d-1}\leq_{lex}x_j w\,''/x_n=x_jx_n^{d-2}\leq_{lex}x_2x_n^{d-2}\leq_{lex} u\,'.
$
Therefore $x_jw\,''\in I\,'$ for $2\leq j\leq n$, hence 
$I\,':w\,''=(x_2,\ldots,x_n).
$
In conclusion we have proved that in every case one may find a monomial $w\notin I\,'$ such that $I':w=(x_2,\ldots, x_n)$.

{\textit{Subcase $C_4$}}: Finally, let $v\leq_{lex} x_2^{d-1}x_n$. In this case, the ideal $I\,':x_2^{d-1}$ obviously contains $(x_2,\ldots,x_n)$. Since the other inclusion is trivial, we get $I\,':x_2^{d-1}=(x_2,\ldots,x_n)$. It is clear that $x_2^{d-1}\notin I\,'$.
\end{proof}

By using Lemma \ref{depthneq0} we get:

\begin{Proposition}\label{depth} Let $I=(\mathcal{L}(u,v))$ be a lexsegment ideal defined by the monomials $u=x_1x_l^{a_l}\ldots x_n^{a_n},\ v=x_q^{b_q}\ldots x_n^{b_n}$ where $a_l,\ b_q>0$, $l,q\geq2$ and $x_n u/x_1<_{lex} v$. Then the following statements hold:
\begin{itemize}
	\item[(a)] If $v=x_2^d$ and $l\geq4$ then $\depth(S/I)=l-2$; 
	\item[(b)] If $v=x_2^{d-1}x_j$ for some $3\leq j\leq n-2$ and $l\geq j+2$ then $\depth(S/I)=l-j$;
	\item[(c)] $\depth(S/I)=1$ in all the other cases.
\end{itemize}
\end{Proposition}

\begin{proof} Since $x_1-x_n$ is regular on $S/I$ if $x_nu/x_1<_{lex}v$, we have $\depth(S/I)=\depth(S\,'/I\,')+1$. The conclusion follows  applying Lemma \ref{depthneq0}.
\end{proof}

\begin{Corollary}
Let $I=(\mathcal{L}(u,v))$ be a lexsegment ideal defined by the monomials $u=x_1x_l^{a_l}\ldots x_n^{a_n},\ v=x_q^{b_q}\ldots x_n^{b_n}$ where $a_l,\ b_q>0$, $l,q\geq2$ and $x_nu/x_1<_{lex} v$. Then the following statements hold:
\begin{itemize}
	\item[(a)] If $v=x_2^d$ and $l\geq4$ then $\projdim(S/I)=n-l+2$; 
	\item[(b)] If $v=x_2^{d-1}x_j$ for some $3\leq j\leq n-2$ and $l\geq j+2$ then $\projdim(S/I)=n-l+j$;
	\item[(c)] $\projdim(S/I)=n-1$ in all the other cases.
\end{itemize}
\end{Corollary}

As a consequence of the results of this section we may characterize the Cohen--Macaulay lexsegment ideals.

In the first place, we note that the only Cohen--Macaulay lexsegment ideal such that $\dim(S/I)=0$ is $ I=\frak m^d$. Therefore it remains to consider Cohen--Macaulay ideals $I$ with $\dim(S/I)\geq1$.

\begin{Theorem} Let $n\geq3$ be an integer, let $u=x_1^{a_1}\ldots x_n^{a_n}$, $v=x_1^{b_1}\ldots x_n^{b_n}$, with $a_1>b_1\geq0,$   monomials of degree $d,$  and $I=(\mathcal{L}(u,v))\subset S$  the lexsegment ideal defined by $u$ and $v.$ We assume  that $\dim(S/I)\geq 1$. Then $I$ is Cohen--Macaulay if and only if one of the following conditions is fulfilled:
\begin{itemize}
	\item[(a)] $u=x_1x_n^{d-1}$ and $v=x_2^d$;
	\item[(b)] $v=x_{n-1}^ax_n^{d-a}$ for some $a>0$ and $x_n\ u/x_1<_{lex}v$. 
\end{itemize}
\end{Theorem}
\begin{proof} Let $u,v$ be as in (a). Then $\dim(S/I)=n-2$, by Proposition \ref{dim} and $\depth(S/I)=n-2$ by using (a) in Proposition \ref{depth} for $n\geq4$ and (c) for $n=3$.

Let $u,\ v$ as in (b). Then $\dim(S/I)=1$ by Proposition \ref{dim}. By using Proposition \ref{depth}(c), we obtain  $\depth(S/I)=1$, hence $S/I$ is Cohen--Macaulay.

For the converse, in the first place, let us take $I$ to be Cohen--Macaulay of $\dim(S/I)=1$. By Proposition \ref{dim} we have $q=n$ or $q=n-1$. If $q=n$, then $v=x_n^d$ and $x_n u/x_1\geq_{lex}v$. By Proposition \ref{depthzero}, $\depth(S/I)=0$, so $I$ is not Cohen--Macaulay.

Let $q=n-1$, that is $v=x_{n-1}^ax_n^{d-a}$ for some $a>0$. By Proposition \ref{depthzero}, since $\depth(S/I)>0$, we must have $x_nu/x_1<_{lex}v$, thus we get (b).

Finally, let $\dim(S/I)\geq2$, that is $q\leq n-2$. By using Proposition \ref{depth}, we obtain  $q=2$. Therefore $\dim(S/I)=\depth(S/I)=n-2$. Using again Proposition \ref{depth} (a),(b), it follows that $u=x_1x_n^{d-1}$ and $v=x_2^d$. 
\end{proof}

\end{document}